\catcode`\@=11

\magnification=1200
\baselineskip=14pt

\pretolerance=500    \tolerance=1000 \brokenpenalty=5000

\catcode`\;=\active
\def;{\relax\ifhmode\ifdim\lastskip>\z@
\unskip\fi\kern.2em\fi\string;}

\overfullrule=0mm

\catcode`\!=\active
\def!{\relax\ifhmode\ifdim\lastskip>\z@
\unskip\fi\kern.2em\fi\string!}

\catcode`\?=\active
\def?{\relax\ifhmode\ifdim\lastskip>\z@
\unskip\fi\kern.2em\fi\string?}

\frenchspacing

\newif\ifpagetitre            \pagetitretrue
\newtoks\hautpagetitre        \hautpagetitre={ }
\newtoks\baspagetitre         \baspagetitre={1}

\newtoks\auteurcourant        \auteurcourant={M. Laurent   }
\newtoks\titrecourant
\titrecourant={On inhomogeneous Diophantine approximation and Hausdorff dimension }

\newtoks\hautpagegauche       \newtoks\hautpagedroite
\hautpagegauche={\hfill\sevenrm\the\auteurcourant\hfill}
\hautpagedroite={\hfill\sevenrm\the\titrecourant\hfill}

\newtoks\baspagegauche       \baspagegauche={\hfill\rm\folio\hfill}

\newtoks\baspagedroite       \baspagedroite={\hfill\rm\folio\hfill}

\headline={
\ifpagetitre\the\hautpagetitre
\global\pagetitrefalse
\else\ifodd\pageno\the\hautpagedroite
\else\the\hautpagegauche\fi\fi}

\footline={\ifpagetitre\the\baspagetitre
\global\pagetitrefalse
\else\ifodd\pageno\the\baspagedroite
\else\the\baspagegauche\fi\fi}

\def\date{\ {\the\day}\
\ifcase\month\or Janvier\or F\'evrier\or Mars\or Avril
\or Mai \or Juin\or Juillet\or Ao\^ut\or Septembre
\or Octobre\or Novembre\or D\'ecembre\fi\
{\the\year}}

\def\up#1{\raise 1ex\hbox{\sevenrm#1}}

\def\cqfd{\unskip\kern 6pt\penalty 500
\raise -2pt\hbox{\vrule\vbox to 10pt{\hrule width 4pt
\vfill\hrule}\vrule}\par\medskip}

\def\section#1{\vskip 7mm plus 20mm minus 1.5mm\penalty-50
\vskip 0mm plus -20mm minus 1.5mm\penalty-50
{\bf\noindent#1}\nobreak\smallskip}

\def\subsection#1{\medskip{\bf#1}\nobreak\smallskip}

\def\displaylinesno #1{\dspl@y\halign{
\hbox to\displaywidth{$\@lign\hfil\displaystyle##\hfil$}&
\llap{$##$}\crcr#1\crcr}}

\def\ldisplaylinesno #1{\dspl@y\halign{
\hbox to\displaywidth{$\@lign\hfil\displaystyle##\hfil$}&
\kern-\displaywidth\rlap{$##$}
\tabskip\displaywidth\crcr#1\crcr}}

\def\hfl#1#2{\smash{\mathop{\hbox to 12 mm{\rightarrowfill}}
\limits^{\scriptstyle#1}_{\scriptstyle#2}}}

%%%%%%%%%%%%%%%%%%%%%%%%%%%%%%%%%%%%%%%%%%%%%%%%%%%%%%%%%%%%%%%%%
\catcode`\@=12

\def\bQ{{\bf Q}}
\def\bR{{\bf R}}
\def\bT{{\bf T}}

\def\bZ{{\bf Z}}

\def\cB{{\cal B}}
\def\cC{{\cal C}}

\def\cH{{\cal H}}
\def\cK{{\cal K}}

\def\del{\delta}

\def\th{{\theta}}

\def\bT{{\bf T}}

\def\bZ{{\bf Z}}

\def\rk{\mathop{\rm rk}\nolimits}

\def\proof{\bigskip\noindent{\it Proof.}\ }

\def\and{\quad\hbox{and}\quad}

\def\om{{\omega}}
\def\omc{{\hat{\omega}}}

%\dimen1=\ht1 \advance\dimen1 by 2pt \dimen2=\dp1 \advance\dimen2 by 2pt
%\setbox1=\hbox{\vrule height\dimen1 depth\dimen2\box1\vrule}%
%\setbox1=\vbox{\hrule\box1}%
%\advance\dimen1 by .4pt \ht1=\dimen1
%\advance\dimen2 by .4pt \dp1=\dimen2 \box1\relax

\def\M{\mathop{\rm M\kern 1pt}\nolimits}
\def\h{\mathop{\rm h\kern 1pt}\nolimits}

\def\romain#1{\uppercase\expandafter{\romannumeral #1}}

\def\card{\mathop{\rm Card\kern 1.3 pt}\nolimits}
\def\deg{\mathop{\rm deg\kern 1pt}\nolimits}
\def\det{\mathop{\rm det\kern 1pt}\nolimits}

\def\h{\mathop{\rm h\kern 1pt}\nolimits} \long\def\forget#1\endforget{}

\def\og{\leavevmode\raise.3ex\hbox{$\scriptscriptstyle 
\langle\!\langle\,$}}
\def\fg{\leavevmode\raise.3ex\hbox{$\scriptscriptstyle
\!\rangle\!\rangle\,\,$}}

%%%%%%%%%%%%%%%%%%%%%%%%%%%%%%%%%%%%%%%%%%%%%%%%%%%%%
%%%%%%%%%%%% Biblio %%%%%%%%%%%%%%%¢

\def\BeDo{1}
\def\BuA{2}
\def\Bu{3}
\def\BuCh{4}
\def\BuLa{5}
\def\BuLaB{6}
\def\Cas{7}
\def\Ch{8}
\def\Fa{9}
\def\Khi{10}
\def\Tr{11}

\centerline{}

\vskip 4mm

\centerline{
\bf  On inhomogeneous Diophantine approximation}
\centerline{
\bf  and Hausdorff dimension}

\vskip 8mm
\centerline{ by Michel L{\sevenrm AURENT}\footnote{}{\rm
} }

\vskip 9mm

\noindent {\bf Abstract --}    Let $\Gamma = \bZ A +\bZ^n \subset \bR^n$ be a dense subgroup with rank $n+1$ and let $\omc(A)$
denote the exponent of uniform simultaneous  rational approximation to the point $A$. We show that for any real number
$v\ge \omc(A)$, the Hausdorff dimension of the set $\cB_v$ of points in $\bR^n$ 
which are $v$-approximable with respect to $\Gamma$,  is equal to $1/v$. 

\vskip15mm

\section
{1. Inhomogeneous approximation.}

We first introduce the general framework of inhomogeneous approximation, following the traditional  setting  employed
in the book of Cassels [\Cas],
and adhering to the notations of [\BuLa] for  the various exponents of approximation involved. 

Let $m$ and $n$ be positive integers and let $A$ be a $n\times m$ matrix with real entries. The transposed matrix of $A$ is 
denoted by ${}^tA$. We consider both the subgroup 
$$
\Gamma = A \bZ^m + \bZ^n \subset \bR^n,
$$
generated modulo $\bZ^n$ by the $m$ columns of $A$,  and its  {\it dual} subgroup $$
\Gamma'  = {}^tA \bZ^n + \bZ^m \subset \bR^m,
$$
generated modulo $\bZ^m$ by the $n$ rows of $A$. It may be enlightening to view alternatively  $\Gamma$ as a subgroup of classes 
modulo $\bZ^n$, lying in the $n$-dimensional torus $\bT^n= (\bR/\bZ)^n$. Kronecker's theorem asserts that $\Gamma$ is dense in
$\bR^n$ iff the dual group $\Gamma'$ has maximal rank $m+n$ over $\bZ$. We shall assume from now that $\rk_\bZ\Gamma'= m+n$.

In order to measure how sharp is the approximation to a given point $\beta$ in $\bR^n$ by elements of $\Gamma$, we introduce 
the following exponent $\om(A,\beta)$. For any point $\th$ in $\bR^n$,  denote by $| \th |$ the supremum  norm of $\th$ and 
by $\| \th \|= \min_{x\in \bZ^n}|\th -x|$ the distance in $\bT^n$ between  $\th$ mod $\bZ^n$ and $0$.

\proclaim 
Definition 1. For any $\beta \in \bR^n$, let  $\om(A, \beta)$ be the supremum, possibly infinite,  of the real numbers $\om$ 
for which there exist infinitely many 
integer points $q\in \bZ^m$ such that 
$$
\| Aq -\beta\| \le | q|^{-\om}.
$$ 

It is plain from the definition that  $\om(A,\beta)\ge 0$. Now, in relation with  the linear independence of the rows of $A$, we 
introduce for any real matrix $M$ the following {\it uniform} homogeneous  exponent:

\proclaim 
Definition 2. Let $M $ be an $m \times n$ matrix with real entries. We denote by  $\omc(M)$  the supremum, possibly infinite,  
of the real numbers $\om$ such that  for any sufficiently large positive real number $Q$,  there exists a non-zero
integer point $q\in \bZ^n$ such that 
$$
   | q | \le Q \and        \| Mq \| \le Q^{-\om}.
$$ 

Dirichlet's box principle shows that $\omc(M) \ge n/m$. We are now able to  formulate the  classical transfer between homogeneous and inhomogeneous
approximation   in terms of these exponents thanks to  the

\proclaim
Theorem 1 \hbox{\rm [\BuLa]}.  For any $n$-tuple $\beta$ of real numbers,  the lower
bound
$$
\om(A,\beta) \ge {1\over \omc({}^tA)} \leqno{(1)}
$$
holds true. Moreover  we have equality of both members  in  (1) for almost all $\beta$ with respect to the
      Lebesgue measure on $\bR^n$.

We come now to our main  topic which is the study for $v\ge 0$ of the family of subsets
$$
\cB_v = \left\{ \beta \in \bR^n ; \quad \om(A, \beta) \ge v \right\} \subseteq \bR^n, 
$$
and of their Hausdorff dimension $\del(v)$ as a function of $v$.  It follows immediately from Theorem 1 that $\cB_v = \bR^n$
when $v\le 1/\omc({}^tA)$, while $\cB_v$ is a null set for $v > 1/\omc({}^tA)$. Furthermore, we know that these latter sets  are
rather small thanks to the following crude result,  quoted as Proposition 7 in [\BuLa]:

\proclaim Theorem 2. For any real   number $v > 1/\omc({}^tA)$, the   Hausdorff dimension $\del(v)$
is strictly less than ~$n$.

In fact, the proof of Proposition 7 of [\BuLa] gives  the explicit upper bound
$$
\del(v) \le n-1  + { 1 \over 1 + ( v \, \omc({}^tA) -1)/ (1+v)}. \leqno{(2)}
$$
On the other hand, an easy application of Hausdorff-Cantelli Lemma (see [\BeDo, \Bu]) provides us  with the following   bound:

\proclaim Theorem 3. For any $v > 0$, we have 
$$
\del(v) \le \min\left( n, {m\over v}\right) . \leqno{(3)}
$$

We refer to Theorem 5 of [\BuCh] for a proof of the inequality (3). 
Note that (2) is certainly sharper than (3) when $v$ belongs to the
interval $[ 1/\omc({}^tA), m/n]$, while the upper bound  (3) is expected to be an equality for sufficiently large values of $v$.
When $m=n=1$, it has been proved independently in [\BuA] and in [\Tr] that $\del(v) = \min(1,1/v)$, so that (3) is indeed an equality
for any $v>0$ in that case. However, the examples displayed in Theorem 1  of [\BuCh] for $(m,n) = (2,1)$ or $(m,n)= (3,1)$, 
 show that the inequality (3) may be strict for any given $v>1$.  
Motivated by Theorem 5  below, we address  the following 

\proclaim Problem.  
Assume that $\omc(A)$ is finite. Show that  $\del(v) = m/v$ for any $v$ sufficiently large in term of $\omc(A)$.  

Notice  that $\omc(A) \ge m/n$. It seems plausible  that the assumption  $v\ge \omc(A)$ should always be sufficient
in order to ensure that  $\del(v) = m/v$. It holds true when $m=1$
according to Theorem 5 below. Note also that the lower bound $v\ge \omc(A)$
 occurs  naturally in  the construction of a Cantor-type  set $\cK$ as in Section 4.

\section
{2. Simultaneous approximation.}

Our knowledge concerning the Hausdorff dimension $\del(v)$ is more substantial  for  $m=1$, that is to say when
$$
\Gamma = \bZ  \left(\matrix{\alpha_1\cr\vdots\cr\alpha_n\cr}\right) + \bZ^n
$$
is generated by a single  vector spinning in $\bT^n$, thanks to the fine results [\BuCh] obtained by Bugeaud and Chevallier. 
 With regard to the above Problem, let us first quote  their Theorem 3 as follows: 
  
\proclaim Theorem 4.  Let $A= {}^t(\alpha_1, \dots , \alpha_n)$ be an $n \times 1$ real matrix with
$1, \alpha_1, \dots, \alpha_n$ linearly independent over $\bQ$. Then $\del(v) = 1/v$ for any $v\ge 1$. 

We state now  our main result.

\proclaim Theorem 5. 
Let $A= {}^t(\alpha_1, \dots , \alpha_n)$ be an $n \times 1$ real matrix with
$1, \alpha_1, \dots, \alpha_n$ linearly independent over $\bQ$. Then the equality $\del(v) = 1/v$ 
holds true for any $v \ge \omc(A)$.

Note that Theorem 5  extends the previous statement since 
$$
{1\over n} \le \omc(A) \le 1.
$$
 The lower bound $\omc(A) \ge 1/n$ follows immediately from Dirichlet's box principle, while the
 upper bound $\omc(A) \le 1$ is implicitely contained in the seminal  work [\Khi] of  Khintchine. 
 It  is expected  that any intermediate value should be reached for some $n\times 1$ matrix $A$.
 We direct to  [\BuLa, \BuLaB] for more precise informations on that topic.

 Theorem 5 implies the following 

\proclaim Corollary. Assume that $\omc(A)=1/n$. Then 
$$
\del(v) = \min\left( n , {1\over v}\right)
$$
for any $v >0$.

The above statement   was initially established by  Bugeaud and  Chevallier in [\BuCh],
 under the stronger assumption that $A$ is  a {\it regular} matrix (according to the terminology of [\Cas]), meaning that there exists a positive
real number $\epsilon$ such that the lower bound
$$
\min_{{q\in \bZ\atop  0 <  q  \le Q}} \| q A \|  \ge \epsilon Q^{-1/n}
$$
holds for arbitrary large values of $Q$.

The proof of Theorem 5  is based on the mass distribution principle [\Bu, \Fa]. This method enables us to bound from below
 the Hausdorff measure $\cH^f(\cB_v)$
of the set $\cB_v$ for suitable dimension functions $f$. It turns out that $\cH^f(\cB_v)=+\infty$  when
$f(r) =  r^{1/v} \log(r^{-1})$  and $v> \omc(A)$, 
as it can be easily seen with some minor modifications of the proof given in Section 4.
Since  $v\mapsto 1/v$ is a decreasing function, a standard argument of Hausdorff measure (see [\BeDo] p. 71) then shows that the 
Hausdorff dimension of the smaller subset
$$
\cB'_v = \left\{ \beta \in \bR^n ; \quad \om(A, \beta) = v \right\} \subseteq \bR^n, 
$$
coincides with the Hausdorff dimension $\del(v)= 1/v$ of $\cB_v$ if $v> \omc(A)$. It follows that for fixed $A$,
the set of values of the exponent
$\om(A,\beta)$ contains the whole interval $]\omc(A), +\infty[$, when $\beta$ ranges over $\bR^n$. 

\section
{2. Best  approximations.}

We review here some properties of the {\it best approximations} to $A$ which are needed for proving Theorem 5. Their detailled proof
can be found in Section 5 of [\BuCh] and in [\Ch]. Throughout this section, $A$ stands for a $n\times 1$ matrix.

A best approximation to $A$ is a positive integer $q$ such that $ \| p A\| > \|ÊqA\|$ for every integer $p$ with $0 < p < q$.
Let $(q_k)_{k\ge 0}$ be the ordered sequence of these best approximations, starting with $q_0=1$. Put
$$
\rho_k = \min_{ 0 < q < q_k} \| q A\| = \| q_{k-1}A\| .
$$
It is readily observed that $ \omc(A)$ is equal to  the lower  limit of the ratio $ \log (\rho_k^{-1})/\log q_{k}$,
 as $k\rightarrow + \infty$. Therefore, if $v$ is any given real number greater than $ \omc(A)$, the inequality
$$
\| q_{k-1} A\| \ge 4 q_k^{-v} \leqno{(4)}
$$
holds for infinitely many $k$.

The key point is to remark that, for large $k$,  the set 
$$
\Gamma_k = \{ qA + \bZ^n; \quad 0 \le  q < q_k\},
$$
when viewed as a subset of $\bT^n$,  is closed to a finite group $\Lambda_k$ which is  well distributed in the torus.
Let $P_k$ be the closest integer point to $q_kA$. Set now 
$$
\Lambda_k = \{ q {P_k\over q_k} + \bZ^n ; \quad 0 \le q < q_k\}= \{ q {P_k\over q_k} + \bZ^n ; \quad q\in \bZ \}.
$$
Clearly $\Lambda_k$  is lattice in $\bR^n$ with determinant $q_k^{-1}$.  Let $\lambda_{1,k} \le \cdots \le \lambda_{n,k}$
be the successive minima of the lattice $\Lambda_k$ with respect to the unit ball $| x| \le 1$.  

\proclaim
Lemma 1. For any integer $k$ and any ball $B(x,r)\subset \bR^n$ centered at the point $x$ with radius $r$, we have the following upper bounds
 \footnote{{\rm(\dag)}}{\rm The constants involved in the symbols $\ll$ and $\asymp$ depend only on  $n$. The ball
 $B(x,r)$ denotes  the  hypercube of points $y\in \bR^n$ with $| y-x | \le r$.}.
If $ r \le \lambda_{i,k}$ for some $i\le n$, then  
$$
 \card \Big(\Gamma_k \cap B(x,r)\Big) \, \ll  \, \prod_{j=1}^{i-1}{ r\over \lambda_{j,k}} \,  \ll \,
 \Big(q_k  \prod_{j=i}^{n}\lambda_{j,k}\Big) r^{i-1}\, ,
$$
(with the convention that the empty product is equal to $1$ when $i=1$). If $r \ge \lambda_{n,k}$, then
$$
 \card \Big(\Gamma_k \cap B(x,r) \Big) \, \ll  \, q_k r^n.
 $$
 Furthermore $\rho_k \asymp \lambda_{1,k}$,  and the last minimum $\lambda_{n,k}$ tends to $0$ when $k$ tends to infinity.

\proof We first prove the above inequalities for $x=0$ with $\Gamma_k$ replaced by $\Lambda_k$. To that purpose, thanks to LLL algorithm,
we use a {\it reduced}  basis $\{e_1, \dots , e_n \}$ of the lattice $\Lambda_k$, meaning that $| e_i | \asymp \lambda_{i,k}$ for $1 \le i \le  n$
and $| \sum x_i e_i | \asymp \max |  x_i e_i |$. We easily obtain the expected bounds for $ \card\Big(\Lambda_k \cap B(0,r)\Big)$, 
using morever Minkowski's theorem on successive minima: 
$$ 
\prod_{j=1}^{n} \lambda_{j,k} \asymp \det\Lambda_k=q_k^{-1}.
$$
See [\BuCh] for more details.
Next, the same  inequalities  hold for any point $x\in \bR^n$ since $\Lambda_k$ is a group.
 In order to replace finally $\Lambda_k$ by $\Gamma_k$,
observe that the distance between the points $qA$ and $qP_k/q_k$ is smaller than $  \rho_{k+1}< \rho_k \ll \lambda_{1,k}$, 
for any integer $q$ with $0\le q < q_k$.

As for the assertions concerning $\lambda_{1,k}$ and $\lambda_{n,k}$, we refer to \S5 of [\BuCh].
 \cqfd
 
 \section
 {4. Proof  of Theorem 5  and of its corollary.}
 
 Let us first deduce the corollary from Theorem 5.  Thanks to  transfer inequalities between uniform exponents due to 
 Apfelbeck and Jarn\'\i k (see  for instance formula (6) in [\BuLa]), we know that $\omc(A) = 1/n$ iff $\omc({}^tA)= n$.
 Then, it  follows from Theorem 1 that $\cB_v = \bR^n$ when $ v \le 1/n$, so that   $\del(v)=n$ for any $v$ in the interval
 $[0,1/n]$. On the other hand,  Theorem 5 gives $\del(v) =1/v$ for $v\ge 1/n$. Therefore, the formula 
 $$
\del(v) = \min\left( n , {1\over v}\right)
$$
 holds true for any positive real number $v$.  
 
 As for the proof of Theorem 5,  note that the dimension $\del(v)$ is a non-increasing function of $v$
 and that $\del(v)\le 1/v$ by Theorem 3. Thus, 
    it suffices to establish  the lower bound 
 $\del(v) \ge 1/v$  for any $v > \omc(A)$. We closely follow the lines of [\BuCh].
 
 Let $v$ and $s$ be positive real numbers such that $ v > \omc(A)$ and $ s < 1/v$. We construct a Cantor-type
 set $\cK \subseteq \cB_v$ whose Hausdorff dimension is  $\ge s$. Let $(k_j)_{j\ge 0}$ be an increasing sequence of positive integers
 such that (4) holds  for any integer $k=k_j, j \ge 0$, appearing in the sequence. The sequence $(k_j)$ is also assumed to be very lacunary,
 in the sense that each   value $k_{j+1}$ is taken  sufficiently large in term of the preceding value $k_j$. 
 The precise meaning of these growth conditions  will be explicited in the course of the construction.
 
 The set $\cK$ is the intersection
 $$
 \cK = \cap_{j\ge 0} K_j
 $$
 of  nested  sets $ K_j$. Each $K_j$ is a finite union of closed balls $B$ with radius $q_{k_j}^{-v}$, 
 centered at some point of $\Gamma_{k_j}$. Therefore $\cK$ is clearly contained in $\cB_v$. Note that the  $K_j$ are  made up 
 with disjoint balls, as a consequence of  (4). 
 We start by taking $k_0$ arbitrary and by choosing for $K_0$ a single ball of the required  type. Put  $N_0=1$. We define inductively
 $K_1 \supset K_2 \supset \dots$ as follows. Suppose that $K_j$ has already been constructed. Since the sequence of points 
 $(qA )_{q\ge 1}$ is uniformely distributed modulo $\bZ^n$ in $\bT^n$ ([\Cas] Chapter IV), we may choose $k_{j+1}$ large enough
 so that each ball occurring in $K_j$, whose Euclidean volume is equal to $2^n  q_{k_j}^{-nv}$,
  contains $\sim 2^n q_{k_{j+1}}q_{k_j}^{-nv}$  points of $\Gamma_{k_{j+1}}$.
  Dropping eventually some  of them, we select in each ball $B$ occurring in $K_j$ exactly the same number
  $$
  N_{j+1} = \Big[ 2^{n-1} q_{k_{j+1}} q_{k_j}^{-nv}\Big]
  $$
  of  points in  $B\cap \Gamma_{k_{j+1}}$ for which the balls $B'$ with radius $q_{k_{j+1}}^{-v}$ centered 
  at these points are included in $B$.  We define  $K_{j+1}$ as the union of  all these selected balls $B'$,  for any  $B$ in $K_j$.
  
  We define now a probability measure  $\mu$ on $\bR^n$ in the following way. First, if $B$ is one of the balls 
  which is part of a set $K_j$, we set
  $$
  \mu(B) ={1 \over N_0\times \cdots \times N_j},
  $$
so that $\mu(K_j)=1$. For any borelian subset $E$, put
$$
\mu(E) = \inf_\cC \Big( \sum_{ B\in \cC} \mu(B)\Big),
$$
where the infimum is taken over all coverings $\cC$ of $ E\cap \cK$ by disjoint balls $B$ occurring 
in the sets  $K_j, j\ge 0$. Then $\mu$ is a probability measure on $\bR^n$ whose support is contained
in $\cK$ [\Fa]. 

\proclaim
Lemma 2. For any point $x\in \bR^n$ and any sufficiently small radius $r$, we have the upper bound
$$
\mu(B(x,r)) \ll r^s.
$$

\proof
Let $j$ be the index determined by 
$$
q_{k_{j+1}}^{-v} \le r < q_{k_j}^{-v}.
$$
The set $\cK\cap  B(x,r)$ is certainly covered by the collection of all balls $B$ with radius $q_{k_{j+1}}^{-v}$  involved in $K_{j+1}$
which intersect $B(x,r)$. Therefore 
$$
\mu(B(x,r)) \le   \sum_{ B\cap B(x,r) \not= \emptyset } \mu(B)
\le {1 \over N_0\times \cdots \times N_{j+1}}  \card \Big(\Gamma_{k_{j+1}} \cap B(x,r+ q_{k_{j+1}}^{-v})\Big) .
\leqno{(5)}
$$  
We make use of Lemma 1 to bound the right hand side of (5). 

Suppose first that
$$
r + q_{k_{j+1}}^{-v} \le \lambda_{1, k_{j+1}}.
$$
Then Lemma 1 (with $i=1$) gives
$$
\mu(B(x,r))  r^{-s} \ll {(q_{k_{j+1}}^{-v})^{-s}\over N_0\times \cdots \times N_{j+1}}
 \ll { q_{k_{j}}^{nv} \over N_0\times \cdots \times N_{j}}  q_{k_{j+1}}^{sv-1} \ll 1,
$$
provided $q_{k_{j+1}} \ge  (q_{k_{j}}^{nv} / (N_0\times \cdots \times N_{j}))^{1/(1-sv)}$ (note that the exponent $sv-1$ is negative). 

Suppose now that there exists an integer $i$ with $1 \le i < n$, such that 
$$
 \lambda_{i, k_{j+1}} \le r + q_{k_{j+1}}^{-v} \le \lambda_{i+1, k_{j+1}}.
$$
We distinguish two cases, depending on whether $i<s$ or $ i\ge s$.  If $i< s$, using Lemma 1, we get the same bound
$$
\mu(B(x,r))  r^{-s} \ll { r^{i-s}\over (N_0\times \cdots \times N_{j+1}) (\lambda_{ 1, k_{j+1}} \times \dots \times \lambda_{ i, k_{j+1}})}
 \ll { q_{k_{j}}^{nv} \over N_0\times \cdots \times N_{j}}  q_{k_{j+1}}^{sv-1},
$$
since
$$
\lambda_{ i, k_{j+1}} \ge  \dots  \ge \lambda_{ 1, k_{j+1}}\gg \rho_{k_{j+1}} \gg q_{k_{j+1}}^{-v} \and
 r \ge q_{k_{j+1}}^{-v}.
 $$
 When $i\ge s$, Lemma 1 and (5) give  the bounds
 $$
 \eqalign{
\mu(B(x,r))  r^{-s} & \ll { 1 \over N_0\times \cdots \times N_{j+1}} r^{i-s} q_{k_{j+1}} \prod_{ \ell = i+1}^n \lambda_{\ell, k_{j+1}}
\cr
& \ll { q_{k_j}^{nv} \over (N_0\times \cdots \times N_{j})q_{k_{j+1}} }\lambda_{i+1, k_{j+1}}^{i-s} q_{k_{j+1}} \prod_{ \ell = i+1}^n \lambda_{\ell, k_{j+1}}
\cr 
& \ll 
 { q_{k_j}^{nv} \over N_0\times \cdots \times N_{j}} \lambda_{n, k_{j+1}}^{n-s} \ll 1,
 \cr}
 $$
provided $\lambda_{n, k_{j+1}} \le  (q_{k_j}^{nv} /( N_0\times \cdots \times N_{j}))^{-1/(n-s)}$. But we also know from Lemma 1 that
$\lambda_{n, k_{j+1}}$ is arbitrarily small when $k_{j+1}$ is sufficiently large (note that the exponent $n-s$ is positive since
$s < 1/v < 1/\omc(A) \le n$).

Suppose finally that 
$$
 r + q_{k_{j+1}}^{-v} \ge \lambda_{n, k_{j+1}}.
 $$
 Recalling that $r \le q_{k_{j}}^{-v}$, Lemma 1 gives now 
 $$
 \eqalign{
\mu(B(x,r))  r^{-s} & \ll { 1 \over N_0\times \cdots \times N_{j+1}} r^{n-s} q_{k_{j+1}} 
\ll  { q_{k_j}^{nv} \over N_0\times \cdots \times N_{j}}  (q_{k_j}^{-v})^{n-s} 
\cr
& \ll { q_{k_j}^{sv} \over N_0\times \cdots \times N_{j} } \ll  { q_{k_{j-1}}^{nv} \over N_0\times \cdots \times N_{j-1} }q_{k_j}^{sv-1} \ll 1,
 \cr}
 $$
provided $q_{k_{j}} \ge  (q_{k_{j-1}}^{nv} / (N_0\times \cdots \times N_{j-1}))^{1/(1-sv)}$. \cqfd

\bigskip

By the mass distribution principle, Lemma 2 ensures that the Hausdorff dimension of $\cK$ is greater or  equal to $s$. Since
$\cK \subseteq \cB_v$, it follows that $\del(v) \ge s$. Taking now $s$ arbitrarily close to $1/v$, we obtain the lower bound
$\del(v) \ge 1/v$. The proof of Theorem 5 is now complete.

\vskip 1cm

\centerline{\bf References }

\vskip 7mm

\item{[\BeDo]}
        V. I. Bernik and M. M. Dodson,
Metric Diophantine approximation on Manifolds,
Cambridge Tracts in Math. and Math. Phys., vol. 137, Cambridge
University Press, 1999.

\item{[\BuA]}
Y. Bugeaud,
{\it A note on inhomogeneous Diophantine approximation},
Glasgow Math. J. 45 (2003), 105--110.

\item{[\Bu]}
Y. Bugeaud,
Approximation by algebraic numbers,
Cambridge Tracts in Mathematics 160, Cambridge, 2004.

\item{[\BuCh]}
Y. Bugeaud and N. Chevallier,
{\it On simultaneous inhomogeneous Diophantine approximation}, 
Acta Arithmetica, 132 (2006), 97--123.

\item{[\BuLa]}
Y. Bugeaud and M. Laurent,
{\it On exponents of homogeneous and inhomogeneous Diophantine approximation}, 
 Moscow Mathematical Journal, 5, 4 (2005), 747--766.

\item{[\BuLaB]}
Y. Bugeaud and M. Laurent,
{\it On exponents of Diophantine Approximation}, 
in: Diophantine Geometry proceedings, Scuola Normale Superiore Pisa,
Ser. CRM, vol. 4, 2007, 101--118.

\item{[\Cas]}
        J. W. S. Cassels,
An introduction to Diophantine Approximation,
Cambridge Tracts in Math. and Math. Phys., vol. 99, Cambridge
University Press, 1957.

\item{[\Ch]}
 N. Chevallier,
{\it Meilleures  approximations d'un \'el\'ement du tore $\bT^2$ et g\'eom\'etrie des multiples de cet \'el\'ement}, 
Acta Arithmetica, 78 (1996), 19--35.

\item{[\Fa]}
        K. Falconer,
Fractal Geometry: Mathematical Foundations and Applications,
Wiley, 1990.

\item{[\Khi]}
A. Ya. Khintchine,
{\it \"Uber eine Klasse linearer diophantischer Approximationen},
Rendiconti Circ. Mat. Palermo 50 (1926), 170--195.

\item{[\Tr]}
J. Schmeling and S. Troubetzkoy,
{\it Inhomogeneous Diophantine Approximation and Angular Recurrence for
Polygonal Billiards},
Mat. Sbornik 194 (2003), 295--309.

\bigskip

\noindent   {Michel LAURENT}

\noindent 
{Institut de Math\'ematiques de Luminy}

\noindent 
{C.N.R.S. -  U.M.R. 6206 - case 907}

\noindent       {163, avenue de Luminy}

\noindent 
{13288 MARSEILLE CEDEX 9  (FRANCE)}

\noindent 
{\hbox{\tt laurent@iml.univ-mrs.fr}}

\bye